\documentclass[12pt,twoside,english,a4paper]{article}
\usepackage{etex}
\usepackage{caption}
\usepackage[latin1]{inputenc}
\usepackage[intlimits] {amsmath}   
\usepackage{graphicx}
\usepackage{psfrag}
\usepackage{ifthen}
\usepackage{fancyhdr}
\usepackage{rotating}
\usepackage{multirow}
\usepackage{booktabs}              
\usepackage{amssymb}
\usepackage{amsbsy}
\usepackage{bm}                    
\usepackage{bbm}
\usepackage{babel}
\usepackage{theorem}
\usepackage{upgreek}  
\usepackage{pifont}   
\usepackage[active]{srcltx}   

\usepackage{suetterl}  
\newfont{\suetdbl}{suet14 scaled 2000}  

\usepackage{yfonts}  
\newfont{\gothdbl}{ygoth scaled 2000} 
\newfont{\frakdbl}{yfrak scaled 2000}
\newfont{\swabdbl}{yswab scaled 2000}

\usepackage[normalem]{ulem}  
\usepackage[textsize=footnotesize,textwidth=1.7cm]{todonotes}
\usepackage[]{defjs3}      

\usepackage[authoryear]{natbib}
\bibpunct{[}{]}{,}{a}{}{;}
\fancypagestyle{plain}{%
  \fancyhead{}%
  \fancyfoot[c]{\sffamily\thepage}%
}
\makeatletter                           
\def\cleardoublepage{\clearpage\if@twoside \ifodd\c@page\else
  \hbox{}
  \vspace*{\fill}
  \thispagestyle{empty}
  \newpage
  \if@twocolumn\hbox{}\newpage\fi\fi\fi}
\makeatother
\begin{document}
\unitlength1.0cm
\frenchspacing

\thispagestyle{empty}
\ce{\bf \large
%
%
Canonical structure of the LLG equation 
}

\medskip 
\ce{\bf \large
for exponential updates in micromagnetism
}

\vspace{4mm}
\ce{J\"org Schr\"oder$^\ast $, Maximilian Vorwerk}

\vspace{4mm}
\ce{Institute of Mechanics, University of Duisburg-Essen,}
\ce{Universit\"atsstr. 15, 45141 Essen, Germany}
\ce{\small e-mail: 
j.schroeder@uni-due.de \& maximilian.vorwerk@uni-due.de
}    

\vspace{4mm}
\begin{center}
{\bf \large Abstract}
\bigskip

{\footnotesize
\begin{minipage}{14.5cm}
\noindent
In this contribution we propose an exponential update algorithm for magnetic moments appearing in the framework of micromagnetics and the Landau-Lifshitz-Gilbert (LLG) equation. 
This algorithm can be interpreted as the geometric integration on spheres, that a priori satisfy the unit length constraint of the normalized magnetization vector. 
Even though the geometric structures for this are obvious and some works already use an exponential algorithm, to the best of the authors' knowledge, there is no canonical structure of the LLG equation for the exponential update algorithm in micromagnetism.
Tensor algebraic reformulations of the LLG equation allow the canonical representation of the evolution equation for the magnetization, which serves as the basis for different integrators.
Based on the specific structure of the exponential of skew symmetric matrices an efficient update scheme is derived. 
The excellent performance of the proposed exponential update algorithm is demonstrated in representative examples.  

\end{minipage}
}
\end{center}


\section{\hspace{-5mm}. Introduction}
\vspace{-4mm}

In this paper we discuss the numerical integration of the evolution equation of the normalized magnetic moments $\bbm$ in micromagnetism (\cite{Bro:1963:m}). 
This is described by a Landau-Lifshitz-type equation which consists of a precession and a dissipation related damping part, cf. \cite{LanLif:1935:ott} and \cite{Gil:2004:apt}. 
Algorithmic implementations are necessary to calculate the evolution of magnetic moments for complex problems. 
Some algorithms may suffer from instabilities caused by the non-fulfillment of the constraint $\norm{\bbm}=1$, this is the case, e.g., for the implicit Euler scheme. 
An overview on constraining methods is given in \cite{Pro:2001:cm,Cim:2008:tll,LewNig:2003:gio} and \cite{ReiXuSch:2022:acs}.
A simple approach to ensure the fulfillment of the unit length constraint is to utilize a projection scheme, which has been demonstrated in publications as \cite{KruPro:2006:rdi} or \cite{Pro:2001:cm}.
Here, the constraint is not considered during the iterative solution procedure itself, but at the end of each individual time step in terms of a renormalization of the magnetization vectors in the element vertices. 
A similar method was presented by \cite{SriKeiMie:2015:chi}, where the solution procedure is based on a staggered scheme.
Another straightforward approach is presented by \cite{Lan:2008:act,Pro:2001:cm,WanZah:2013:ars,ZhaZhaPei:2016:afe}. 
Here, the authors penalize all magnetization vectors that deviate from the unit length by adding a penalty term to the system of equations. 
The penalization intensity is controlled by a penalty parameter that is especially for complex and strongly heterogeneous microstructures not easy to fit.
A slightly different way to restrict the vectors onto the unit sphere was presented by \cite{SzaBudTouFru:acf:2008:acf,OhmYiGutXu:2022:pfm} and \cite{ReiNieSch:2023:ems} in terms of a Lagrange multiplier. 
Both methods, penalization or renormalization, introduce additional (artificial) contributions to the total energy of the systems.
Therefore, the resulting energies may differ from the expected physical energy of the magnetic system.

This problem is circumvented for a priori length preserving methods. 
In this context spherical coordinates have been applied by \cite{SueSchFid:2000:mso,FidSch:2000:mmt,YiXu:2014:acf}.
Another advantage of spherical coordinates is that they reduce the system of equations by one degree of freedom compared to Cartesian coordinates, what leads to faster simulations, cf. \cite{SchFidSchSueForTsi:2003:spm}. 
However, this advantage is outweighed by the singularities in the poles. 
Magnetization vectors that reach these singularities usually lead to non-converging solution procedures. 
To overcome non-physical switching behavior \cite{DorSchXuKeiMue:2018:cpf}, an extension of the work of \cite{YiXu:2014:acf}, introduced a restriction of the azimuthal angle onto its definition domain $\left(-\pi, \pi\right]$ on the interpolation level.
Also length preserving is the incremental update algorithm for the magnetization vectors proposed by \cite{MieEth:2014:agc}, using the geometrical structure for the integration on spheres which results after several algebraic manipulations in an exponential scheme. 
The authors use the geometrical structure for the integration on spheres and construct, after some clever rearrangements, an exponential update. 
However, they do not rely on a canonical structure for the exponential scheme. 
Further, important aspects, especially concerning the treatment within the framework of the Finite-Element-Method can be found in \cite{DorWoeWul:2023:cbn}.
Another approach to preserve the length of the magnetization vectors was proposed by \cite{Alo:2008:anf}. Here, the authors introduce a so-called $\theta$-scheme in analogy to the heat equation which is unconditional stable for a suitable choice of $\theta$.
A clever reformulation of the initial equation, e.g. based on the mid-point rule, leads to no further numerical methods to constrain the length, as demonstrated in \cite{dAqSerMia:2005:gio}.

Our contribution is mainly motivated by the works of \cite{FenVis:2002:saa,VisFen:2002:qba} and \cite{LewNig:2003:gio}.
In \cite{FenVis:2002:saa} the stability and accuracy of Euler and quaternion micromagnetic algorithms including Langevin random-fields has been investigated in detail. They conclude that their proposed quaternion based scheme is significantly more stable compared to Euler's method.
A quaternion-based algorithm with separate precession rotations and dissipation damping, allowing analytical calculation of the precession part over long time intervals, was proposed by \cite{VisFen:2002:qba}. 
This procedure allows the usage of large time increments compared to conventional algorithms, especially in the case of low anisotropy and weak exchange coupling. 
The philosophy that numerical algorithms should capture the intrinsic (geometric) properties of the system has been analyzed in the well understandable but more mathematically oriented work by \cite{LewNig:2003:gio}.
Starting on the motivation of the usage of geometric integrators, different integrators are discussed and analyzed for micromagnetic applications. 

This work is structured as follows:
\textit{First,} the LLG equation is reformulated and brought into a form suited for our proposed exponential update framework; the canonical form. 
\textit{Second,} performing the algorithmic derivation for the implementation of the exponential method.
\textit{Third,} comparison of the performance of the proposed exponential method with the implicit Euler method and the midpoint rule using representative numerical examples.
\textit{Finally,} the results as well as the advantages and disadvantages of the method are concluded.


\section{\hspace{-5mm}. Canonical structure of the Landau-Lifshitz-Gilbert equation}
\vspace{-4mm}
In dynamic micromagnetism, both an evolution equation and an enthalpy functional are needed to describe the magnetization vectors in space and time.
The evolution equation, here the Landau-Lifshitz-Gilbert equation (LLG), defines how the magnetization vectors evolve under external and internal influences. 
These external and internal influences are described by the magnetic enthalpy and enter the LLG equation in the form of the so-called effective field $\bH^{\rm eff}$, which is derived from a variational derivative of the enthalpy regarding the magnetization. 
In this section, the LLG equation is introduced and reformulated for the discretization discussed in this work. 
Subsequently, the enthalpy functional and the effective field are specified.
For an extensive overview on the foundations on magnetism and micromagnetics in particular the reader is referred to standard literature as \cite{LanLif:1935:ott}, \cite{Bro:1963:m}, \cite{CulGra:2009:itm}, \cite{Coe:2010:mam}, \cite{HubSch:2008:mdt}, \cite{KroPar:2007:hom} and \cite{Ber:1998:him}.

In order to describe the evolution of the magnetic moments
$ \bM = M_s\,\bbm $
, a representation using the magnetic unit director $\bbm$ with the constraint condition 
$\norm{\bbm}=1 $ and $\partial_t \norm{\bbm} =0$ is applied. 
Hence, the so-called semi-implicit form of the LLG equation appears as
\begin{equation}
\partial_t\bbm = 
	-\gamma\, \bbm \times \bH^{eff} 
	    + \alpha \, \bbm \times\partial_{t}\bbm
\label{eq:LLGsemiimplicit}
\end{equation}
with 
$ \gamma = \mu_0 \, \gamma_e 
   \approx 2.2128 \cdot 10^{-5} \, \textrm{m/As} $.
The goal is to derive a canonical form for the evolution equation for the magnetization
\eb
\partial_t {\bbm} = \bW (\bbm, \bH) \cdot \bbm  \, ,
\ee 
suitable for exponential update algorithms.
Substituting the time derivative of the magnetization appearing 
on the right-hand side of (\ref{eq:LLGsemiimplicit})
by the right-hand side of (\ref{eq:LLGsemiimplicit}) yields 
\begin{equation}
\partial_t\bbm = 
	-\gamma\, \bbm \times \bH^{eff} 
	    - \alpha \, \gamma\, \bbm \times
	    ( \bbm \times \bH^{eff} )
	     + \alpha^2 \, \bbm \times (\bbm \times\partial_{t}\bbm) \, .
\label{eq:LLGexplicita}
\end{equation}
Applying the Gra\ss mann identity to the last term of the right-hand-side of \ref{eq:LLGexplicita} and taking into account 
$\partial_t \norm{\bbm}^2 = 2 \, \partial_t {\bbm} \cdot \bbm =  0$,
 yields 
\eb
\bbm \times (\bbm \times \partial_{t}\bbm) = - \partial_{t}\bbm \, .
\ee
Thus the evolution equation for the magnetization appears  
in the so-called explicit form of the LLG equation
\eb
\partial_t {\bbm} = 
-\dfrac{\gamma}{1 + \alpha^2} \bbm \times \bH^{eff}
-\dfrac{\alpha \, \gamma}{1 + \alpha^2} 
 \bbm \times ( \bbm \times \bH^{eff} ) \, .
\label{eq:LLGexplicit}
\ee
For convenience, we introduce the abbreviation $ \beta = \dfrac{\gamma}{1 + \alpha^2}$
and change the order in the terms of the cross product, i.e.,
\eb 
\partial_t {\bbm} = 
\beta \, \bH^{eff} \times  \bbm 
-  \alpha\beta \, ( \bH^{eff} \times  \bbm) \times \bbm  \, .
\label{eq:LLGexplicitnewW}
\ee
By introducing the skew-symmetric matrix 
\ebn
	\BOmega_{\left[ \bH \right]} = - \BOmega_{\left[ \bH \right]}^T =
	\begin{bmatrix}
		0 & -H_3 & H_2\\
		H_3 & 0 & -H_1\\
		-H_2 & H_1 & 0
	\end{bmatrix} \, 
\een
we substitute 
\eb
\bH^{eff} \times  \bbm  = \BOmega_{\left[ \bH^{eff} \right]} \cdot  \bbm \, .
\ee 
Furthermore, for 
$ 
\BOmega_{\left[ \bH^{} \times  \bbm \right]} = 
\bbm \otimes  \bH - \bH^{} \otimes  \bbm
= 2 \, \textrm{skew} [\BOmega_{\left[ \bH \right]} \cdot \BOmega_{\left[\bbm \right]}]
$, 
with the explicit expression
\ebn
\BOmega_{\left[ \bH^{} \times  \bbm \right]} =
\begin{bmatrix}
 0 &  H_2 m_1 - H_1 m_2 &  H_3 m_1 - H_1 m_3 \\
  -H_2 m_1 + H_1 m_2 & 0 &  H_3 m_2 - H_2 m_3 \\
 -H_3 m_1 + H_1 m_3  & -H_3 m_2 + H_2 m_3 & 0 
	\end{bmatrix} 
 \, ,
\een
we substitute 
\eb
( \bH^{eff} \times  \bbm) \times \bbm  = \BOmega_{\left[ \bH^{eff} \times  \bbm \right]} \cdot \bbm \, .
\ee 
Inserting this relations we obtain form (\ref{eq:LLGexplicitnewW}) the evolution 
equation for the magnetization vector
\eb\boxed{
\partial_t {\bbm} = \bW\cdot  \bbm
\quad\textrm{with}\quad
\bW := \beta \, ( \BOmega_{\left[ \bH^{eff} \right]} 
-  \alpha \,  \BOmega_{\left[ \bH^{eff} \times  \bbm \right]})  ,
\label{eq:LLGexplicitnewabbr}
}
\ee
with the skew-symmetric matrix 
\ebn
\bW = \beta  
\begin{bmatrix}
 0 & -H_3 - \alpha (H_2 m_1 - H_1 m_2) &  H_2 + \alpha ( H_1 m_3 - H_3 m_1) \\
 H_3 + \alpha (H_2 m_1 - H_1 m_2) & 0 &  - H_1- \alpha (H_3 m_2 - H_2 m_3) \\
 -H_2 - \alpha (H_1 m_3 -H_3 m_1)  & H_1 + \alpha (H_3 m_2 - H_2 m_3 ) & 0 
	\end{bmatrix} .
\een
For convenience in further mathematical transformations, 
we introduce the following abbreviations for the terms in $\bW$:
\eb
\bW = 
\begin{bmatrix}
 0       & - W_{21} &   W_{13} \\
W_{21}   &    0     & - W_{32} \\
- W_{13} &  W_{32} & 0 
	\end{bmatrix} 
	 =
	\begin{bmatrix}
		0 & -w_3 & w_2\\
		w_3 & 0 & -w_1\\
		-w_2 & w_1 & 0
	\end{bmatrix} \, 
= \BOmega_{\left[ \bw \right]} \, .
\label{eq:skew-symmetricmatrix}
\ee

\section{\hspace{-5mm}. Integration schemes}

The canonical form of the LLG equation  (\ref{eq:LLGexplicitnewabbr}) is the basis for all further time integrators. 
First, the exponential update algorithm is derived, which is based on the series expansion of a skew-symmetric matrix $\bW $. 
With the help of the explicit representation of the cofactor  of $\bW $ and the application of the Caley-Hamilton theorem, the implicit algorithm to be solved iteratively finally follows. 
The geometric properties of the numerical integrator lead to the preservation of structure (invariants, constraints) at the level of the time-discrete evolution equation. 
For a  general discussion of structure-preserving algorithms we refer 
to \cite{HaiLabWan:2006:gni}.
Subsequently, the implicit Euler algorithm (backward Euler) and the midpoint rule are specified for a comparative study.

\subsection{Exponential update algorithm} 
An appropriate integrator for the evolution equation of the magnetization (\ref{eq:LLGexplicitnewabbr}) in the special orthogonal group $\textrm{SO}(3)$, 
defined by 
$\textrm{SO}(3) := \{ \bQ : \IR^3 \rightarrow  \IR^3  \mid \bQ^T \bQ = \bone \, , \; \det \bQ = 1\}$, is 
\eb
{\bbm}_{n+1} = \textrm{exp}[ \bW^{\Delta t}_{n+\xi} ]\cdot  \bbm_n 
\quad\textrm{with}\quad
\xi \in [0, 1] \, .
\label{eq:expLLGexplicitnewalpha}
\ee
 $\xi = 0$ represents an explicit exponential Euler scheme, 
$\xi = 1/2$ the averaged exponential midpoint rule and
$\xi = 1$ the implicit exponential Euler scheme. 
Without further restrictions, we set $\xi = 1$ for the following considerations, thus we focus on 
\eb
{\bbm}_{n+1} = \textrm{exp}[ \bW^{\Delta t}_{n+1} ]\cdot  \bbm_n
\label{eq:expLLGexplicitnew}
\ee
with the time increment $\Delta t = t_{n+1} - t_n$ and the time discretized matrix 
\eb
\bW^{\Delta t}_{n+1} := \beta \, \Delta t \, ( \BOmega_{\left[ \bH_{n+1} \right]} 
- \alpha \,  \BOmega_{\left[ \bH^{eff}_{n+1} \times  \bbm_{n+1} \right]} ) \, .
\label{eq:wdelt}
\ee
The series expansion for the exponential matrix form, with an infinite number of terms, is given by 
\eb
\textrm{exp}[ \bW^{\Delta t}_{n+1} ] = \bone + \bW^{\Delta t}_{n+1} 
+ \dfrac{1}{2!} (\bW^{\Delta t}_{n+1})^2 
+ \dfrac{1}{3!} (\bW^{\Delta t}_{n+1})^3 + \cdots + \dfrac{1}{n!} (\bW^{\Delta t}_{n+1})^n + ... \, .
\ee
Using the Cayley-Hamilton theorem, i.e.,
\eb
\bW^3 = I_3 \, \bone - I_2 \, \bW + I_1 \bW^2  \, ,
\ee
we recognize that every power of the quadratic matrix can be represented as a tensor matrix with the tensor generators up to the second power:
\eb
\bW^n = \gamma_0^{(n)}(I_1, I_2, I_3) \, \bone 
+ \gamma_1^{(n)}(I_1, I_2, I_3) \, \bW + \gamma_2^{(n)}(I_1, I_2, I_3)\bW^2  \, .
\ee
The parameters $\gamma_i^{(n)}(I_1, I_2, I_3) | i = 1,2,3$ are functions of the 
invariants 
\eb
I_1 = \trb{\bW^{\Delta t}} \, , \; 
I_2 = \dfrac{1}{2} (I_1^2 -\trb{(\bW^{\Delta t})^2}) = \trb{\Cofb{\bW^{\Delta t}}}
\, , 
\; I_3 = \detb{\bW^{\Delta t}} \, 
\ee
with the cofactor of (\ref{eq:skew-symmetricmatrix}) given by 
\ebn
\Cofb{\bW} =
\begin{bmatrix}
W_{32}^2 & W_{13} W_{32} &  W_{21} W_{32}\\
W_{13} W_{32}    & W_{13}^2 &   W_{13} W_{21}\\
W_{21} W_{32}    & W_{13} W_{21} &   W_{21}^2
	\end{bmatrix} 
	 =
\begin{bmatrix}
w_{1}^2 & w_{1}\, w_{2} &  w_{1} \, w_{3}\\
w_{2} \, w_{1}    & w_{2}^2 &   w_{2} \, w_{3}\\
w_{3} \, w_{1}    & w_{3} \, w_{2} &   w_{3}^2
	\end{bmatrix} 
= \bw \otimes \bw \, .
\een
In our specific case, a skew-symmetric matrix, results in the following invariants 
\eb
I_1 = 0\, , \quad 
I_2 = \trb{\Cof{\bW^{\Delta t}}} = \bw \cdot \bw =: w^2
\, , \quad 
I_3 = 0 \, 
\ee
and the simplified relation using the Cayley-Hamilton theorem, i.e.,
\eb
\bW^3 =   - I_2 \, \bW   \, .
\ee
If we now use a recursive function to calculate the odd and even powers of the matrix $\bW$
\eb
\bW^{2n+1} =   (- w^2)^{n} \, \bW   \; ,
\qquad 
\bW^{2n} =   (- w^2)^{n-1} \, \bW^2 \, ,
\label{eq:potenzenW}
\ee
applies for $n \ge 1 $ respectively. 
Thus, $ \textrm{exp}[ \bW^{\Delta t} ]$ appears in series expansion form as
\eb
\textrm{exp}[ \bW^{\Delta t}] = 
  \bone 
+
\sum_{n=0}^{\infty} \dfrac{1}{(2n + 1)!} \; (\bW^{\Delta t})^{2n+1}
+
\sum_{n=1}^{\infty} \dfrac{1}{(2n)!} \; (\bW^{\Delta t})^{2n} \, .
\ee
By inserting (\ref{eq:potenzenW}) into the latter expression, the 
tensor representation simplifies to
\eb
\textrm{exp}[ \bW^{\Delta t}] = 
  \bone 
+
\sum_{n=0}^{\infty} \dfrac{(- w^2)^{n}}{(2n + 1)!} \; \bW^{\Delta t} 
+
\sum_{n=1}^{\infty} \dfrac{(- w^2)^{n-1}}{(2n)!} \; (\bW^{\Delta t})^{2} \, .
\label{eq:alogoexprecurs}
\ee
Let us express (\ref{eq:alogoexprecurs}) in terms of the general tensor function for the exponential function of the skew-symmetric matrix
\eb
\textrm{exp}[ \bW^{\Delta t} ] = 
 \bone +\alpha_1 \bW^{\Delta t} + \alpha_2 (\bW^{\Delta t})^2 \, ,
\label{eq:algoexpalpha}
\ee
with vanishing first and third invariant, we identify
\eab
\begin{array}{rl}
\alpha_1 &= 1 - 
\dfrac{w^2}{3!} + \dfrac{w^4}{5!} - \dfrac{w^6}{7!} + \dfrac{w^8}{9!} - \dfrac{w^{10}}{11!} \pm \cdots \, \quad\rightsquigarrow\quad 
  \alpha_1 = \dfrac{\sin (w)}{w}  \, ,
\\[3mm]
\alpha_2 &= \dfrac{1}{2}   
- \dfrac{w^2}{4!} + \dfrac{w^4}{6!} - \dfrac{w^6}{8!} + \dfrac{w^8}{10!} - \dfrac{w^{10}}{12!}\pm \cdots \, \quad\rightsquigarrow\quad
\alpha_2 = \dfrac{1- \cos (w)}{w^2} \,.
\end{array}
\label{eq:alternative2}
\eae


For the implicit time integrator, the equation for the exponential update must be solved iteratively, so we write the equation in residual form and solve 
\eb
\bR_{exp} := {\bbm}_{n+1} - \textrm{exp}[ \bW^{\Delta t}_{n+1} ]\cdot  \bbm_n
\label{eq:resexpLLGexplicitnew}
\ee
with respect to $\bbm_{n+1}$ until 
$\norm{\bR_{exp}} \le \textrm{tol} $, with the predefined tolerance tol. 
For the numerical treatment we redefine the exponential 
as
\eb
\textrm{exp}[ \bW_{n+1}^{\Delta t} ] = 
 \bone + \alpha_1^\ast \, \bW_{n+1}^{\Delta t} + \alpha_2^\ast \, (\bW_{n+1}^{\Delta t})^2 
\label{eq:algoexpalpha3}
\ee
with case differentiations to avoid singularities
\eb
 \alpha_1^\ast = \dfrac{\sin (w)}{w} 
 \quad\textrm{for}\quad 
 w \ge \epsilon
  \quad\textrm{else}\quad 
  \alpha_1^\ast = 1 - \dfrac{w^2}{3!} + \dfrac{w^4}{5!} - \dfrac{w^6}{7!}
\ee 
and 
\eb
\alpha_2^\ast = \dfrac{1- \cos (w)}{w^2}
 \quad\textrm{for}\quad 
 w^2 \ge \epsilon
  \quad\textrm{else}\quad 
  \alpha_2^\ast = \dfrac{1}{2}   - \dfrac{w^2}{4!} + \dfrac{w^4}{6!} \, ,
\ee
with a predefined small positive values $\epsilon$. 
An excellent description of an alternative approach to calculate the exponential of a $3 \times 3$ matrix based on a recursion formula can be found, e.g., in \cite{Eid:2004:aim}. 

\subsection{Backward Euler and Midpoint rule} 
To evaluate the performance of the presented exponential algorithm, 
we compare the results with the implicit Backward Euler method and the midpoint rule. 
Both algorithms are briefly presented here.
Based on the proposed representation of the LLG equation, the backward Euler algorithm appears as
\eb
{\bbm}_{n+1} \approx \bbm_n + \Delta t \, \dot{\bbm}_{n+1}
\quad\textrm{with}\quad
\dot{\bbm}_{n+1} = \bW_{n+1} \cdot {\bbm}_{n+1} \; .
\label{eq:BE1}
\ee
This leads to the equation
\eb
{\bbm}_{n+1} \approx  \bbm_n + \Delta t \, \bW_{n+1}^{\Delta t} \cdot {\bbm}_{n+1} \, ,
\label{eq:BE1a}
\ee
where we use the abbreviation from (\ref{eq:wdelt}).
Equation (\ref{eq:BE1a}) has to be solved iteratively until the residuum 
\eb
\bR_{BE} := (\bone - \bW^{\Delta t}_{n+1}) \cdot {\bbm}_{n+1}  - \bbm_n 
\label{eq:BE1b}
\ee
is approximately zero, i.e., $\norm{\bR_{BE}} \le \textrm{tol} $.
The midpoint rule is an implicit discrete scheme with verified unconditional convergence, see \cite{BarPro:2007:cpi}.
If we apply the proposed algorithm to our reformulation of the LLG equation, we obtain 
\eb
{\bbm}_{n+1} \approx \bbm_n + \Delta t \, \dot{\bbm}_{n+1}
\quad\textrm{with}\quad
\dot{\bbm}_{n+1} = \bW_{n+1/2} \cdot  {\bbm}_{n+1} \, ,
\label{eq:midpr}
\ee
where we have used the abbreviation 
\eb
\bW_{n+1/2} := \bW [\bH^{eff}_{n+1/2}, \bbm^{}_{n+1/2}]
\quad\textrm{with}\quad
(\bullet)_{n+1/2} = \dfrac{1}{2} ((\bullet)_{n+1} + (\bullet)_{n}) \, .
\label{eq:midpra}
\ee
Based on our notation, we introduce  $\bW^{\Delta t}_{n+1/2} = {\Delta t} \; \bW_{n+1/2}$
and define the residuum  
\eb
\bR_{MR} := (\bone - \bW^{\Delta t}_{n+1/2}) \cdot {\bbm}_{n+1}  - \bbm_n \, ,
\label{eq:BE1b}
\ee
 which must be solved iteratively until $\norm{\bR_{MR}} \le \textrm{tol} $ is satisfied.

\sect{\hspace{-5mm} Magnetic enthalpy and effective field}
\vspace{-4mm}

The magnetic effects to be included within the effective field are stored within the magnetic enthalpy function 
\begin{equation}
	\mathcal{H} ( \bbm, \nabla \bbm, \bH ) 
	= \mathcal{H}^{\rm mag} ( \bbm, \bH )
  + \mathcal{H}^{\rm exc} ( \nabla \bbm )
  +  \mathcal{H}^{\rm ani} ( \bbm ) \, ,
\end{equation}
consisting of 
the sum of the contributions for 
the magnetostatic part $\mathcal{H}^{\rm mag}$,
the exchange energy $\mathcal{H}^{\rm exc}$ and 
the magneto-crystalline anisotropy $\mathcal{H}^{\rm ani}$.
Based on this functional the effective field can be derived as 
\begin{equation}
\bH^{\rm eff}:= \frac{1}{\mu_0} \left( \div \BPi - \partial_{\bM} \mathcal{H} \right), \quad \textrm{with}\quad \BPi = \frac{\partial\, \mathcal{H}}{\partial\,\nabla\bbm} \, .
	\label{eq:effective_field}
\end{equation}
The magnetostatic energy as well as the exchange energy are stated as
\begin{equation}
\mathcal{H}^{\rm mag} ( \bbm, \bH ) 
=  - \frac{1}{2}\mu_0 M_s \, \bbm\cdot \bH- \mu_0 \bH\cdot \bH 
\quad\textrm{and}\quad 
\mathcal{H}^{\rm exc} ( \nabla \bbm ) 
    = A^{\rm exc} \, \abs{\nabla \bbm}^2 \, ,
\end{equation}
with $A^{\rm exc}$ denoting the exchange constant.
For different types of crystalline anisotropies we consider several model problems for
$\mathcal{H}^{\rm ani}$ with $\rm  ani = \{iso, ti, cub \}$.
The expression for an isotropic $\rm ani = iso$, a transversely isotropic  $\rm ani = ti$
and a cubic anisotropy $\rm ani = cub $ potential are given by 
\begin{align}
\begin{array}{l l}
\mathcal{H}^{\rm iso} &= \dfrac{K_1^{\rm iso}}{2}\;\bbm\cdot\bbm, 
\quad
\mathcal{H}^{\rm ti} = 
              \dfrac{K_1^{\rm ti}}{2} \left[1- J_1^2 \right]  
\quad\textrm{and}\quad\\[5mm]
\mathcal{H}^{\rm cub} &= 
K_1^{\rm cub} \left[
J_1^2 \, J_2^2 + J_2^2 \, J_3^2 +J_3^2 \, J_1^2 \right] 
+ 
K_2^{\rm cub} J_1^2 \,  J_2^2 \, J_3^2 \, ,
\end{array}
\label{eq:enthalpy_function}
\end{align}
\begin{equation}
	J_i = \bbm\cdot\ba_i
	\quad\textrm{for}\quad i = 1,2,3
\end{equation}
where the unit vectors $\ba_1=\left[1,0,0\right]^T$, $\ba_2=\left[0,1,0\right]^T$ and $\ba_3=\left[0,0,1\right]^T$ represent the directions of the anisotropy axis.
A typical anisotropy surface for typical cubic magnetocrystalline anisotropy constants 
$K_1, K_2$ is depicted in Fig. \ref{nodal_evolution}.

\begin{Figure}[ht]
\begin{picture}(10,4.)
\unitlength1cm
\put( 3.0,-3.5){\includegraphics[width=0.6\textwidth]{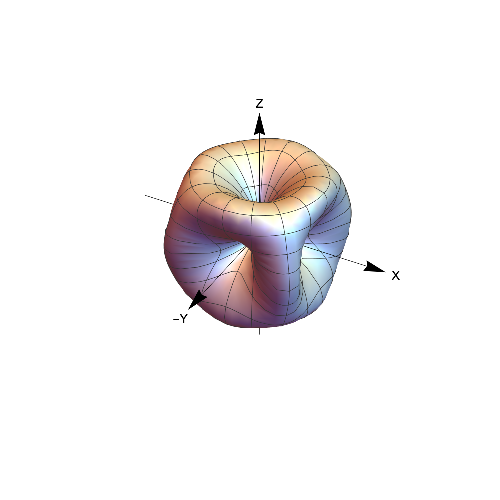}}
\end{picture}
\vspace*{5mm} 
\caption{
Magneto-crystalline energy surface for the cubic anisotropy case ($K_1^{cub} = 2\cdot10^4$ J/m$^3$ and $K_2^{cub} = -4.5\cdot10^4$ J/m$^3$). Parameters taken from \cite{YiXu:2014:acf}.
\label{nodal_evolution}
} 
\end{Figure}

Based on this assumptions the effective field appears as 
\begin{equation}
	\bH^{\rm eff} 
	=  \frac{1}{\mu_0} \left[ 
	\div \left[\frac{\partial \mathcal{H}}{\partial \nabla \bbm} \right]
	- 
	\frac{\partial \mathcal{H}}{\partial \bbm} 
	\right]
=  \underbrace{\frac{ 2 A^{\rm exc}}{\mu_0 M_s} \nabla^2 \bbm}_{\bH^{\rm exc}} 
+ \underbrace{\frac{1}{\mu_0 M_s} \bH}_{\bH^{\rm mag}} 
- \underbrace{\frac{1}{\mu_0} \partial_{\bbm} \mathcal{H}^{\rm ani}}_{\bH^{\rm ani}}
\end{equation}
where $\nabla^2$ denotes the Laplace operator.

\sect{\hspace{-5mm} Comparative numerical studies}
\vspace{-4mm}

The workability and robustness of the developed algorithm is proofed by applying suitable numerical examples for a point element. 
Hence, the evolution of a single magnetization vector $\bbm$ in space and time is considered under different loading scenarios and eventually compared to the performance of the Backward Euler method. 
Within this contribution the external field is considered to be $\bH^{\rm ext}=10^6$A/m.
The material parameters presented in Tab.~\ref{tab:matparam_galf} mimic the soft magnetic material Fe$_{81.3}$Ga$_{18.7}$, also used in \cite{YiXu:2014:acf}.

\begin{Table}[htb!]
	\caption{Material parameters of Fe$_{81.3}$Ga$_{18.7}$ taken from \cite{YiXu:2014:acf}.}
\begin{tabular}{ll|r || ll|r }
  \hline
  Parameter   && Fe$_{81.3}$Ga$_{18.7}$ &Parameter && Fe$_{81.3}$Ga$_{18.7}$ \\[1mm]
  \hline
  \rule{0pt}{13pt}
 anisotropy const.    &$K_1$              $\rm \frac{J}{m^3}$  &  $2 \times 10^4$        &   
 anisotropy const.    &$K_2 $             $\rm\frac{J}{m^3}$   &  $-4.5\times 10^4$       \\[3mm] 
 exchange const.      &$A_{\rm exc}$      $\rm\frac{J}{m}$     &  $\approx 10^{-11}$    & 
 sat. magnetization   &$M_{\rm s}$        $\rm\frac{A}{m}$     &  $1.432 \times 10^6$    \\[3mm]
 vac. permeability    &$\mu_0$            $\rm\frac{H}{m}$     &  $4 \pi \times 10^{-7}$ &
 gyromagnetic ratio   &$\gamma_0$         $\rm\frac{1}{Ts}$    &  $1.76 \times 10^{11}$  \\[1mm]
  \hline
\end{tabular}
\label{tab:matparam_galf}
\end{Table}

Since the magnetization evolution is only considered locally for one point, the exchange energy cannot be included within the effective field. However, to keep completeness of the material parameter set, also the exchange coefficient is given. 

\subsection{Precessional switching}
The switching in a magnetic particle is mainly determined by two terms in the LLG equation. These terms can be differentiated as the precessional part, which makes the magnetization vectors rotate around the effective field, and the dissipative damping part, which makes the magnetization align itself to the effective field, see \cite{Abe:2019:mas}. 
Combined, this exhibits a damped precession, i.e., the magnetization spirally moves  towards the effective field over time until Brown's equation $\bbm\times\bH^{\rm eff}=\bzero$ is eventually satisfied, see \cite{Ber:1998:him}. 
This damped precession is considered throughout this example. The set up of the boundary value problem considered here is analogous to that described in \cite{VisFen:2002:qba} and \cite{YiXu:2014:acf}.  
The magnetization vector initially points into the x$_1$-axis while the magnetic field is applied parallel to the x$_3$-axis. The evolution of the magnetization vectors is studied for different damping coefficients (brown $\alpha$=0.0, gray $\alpha$=0.01, blue $\alpha$=0.1, red $\alpha$=1.0) resulting in different evolution paths towards equilibrium. 

\begin{Figure}[h]
\begin{picture}(0,4.9)
\unitlength1cm
\put(-1.2, -1.0){\includegraphics[width=0.7\textwidth]{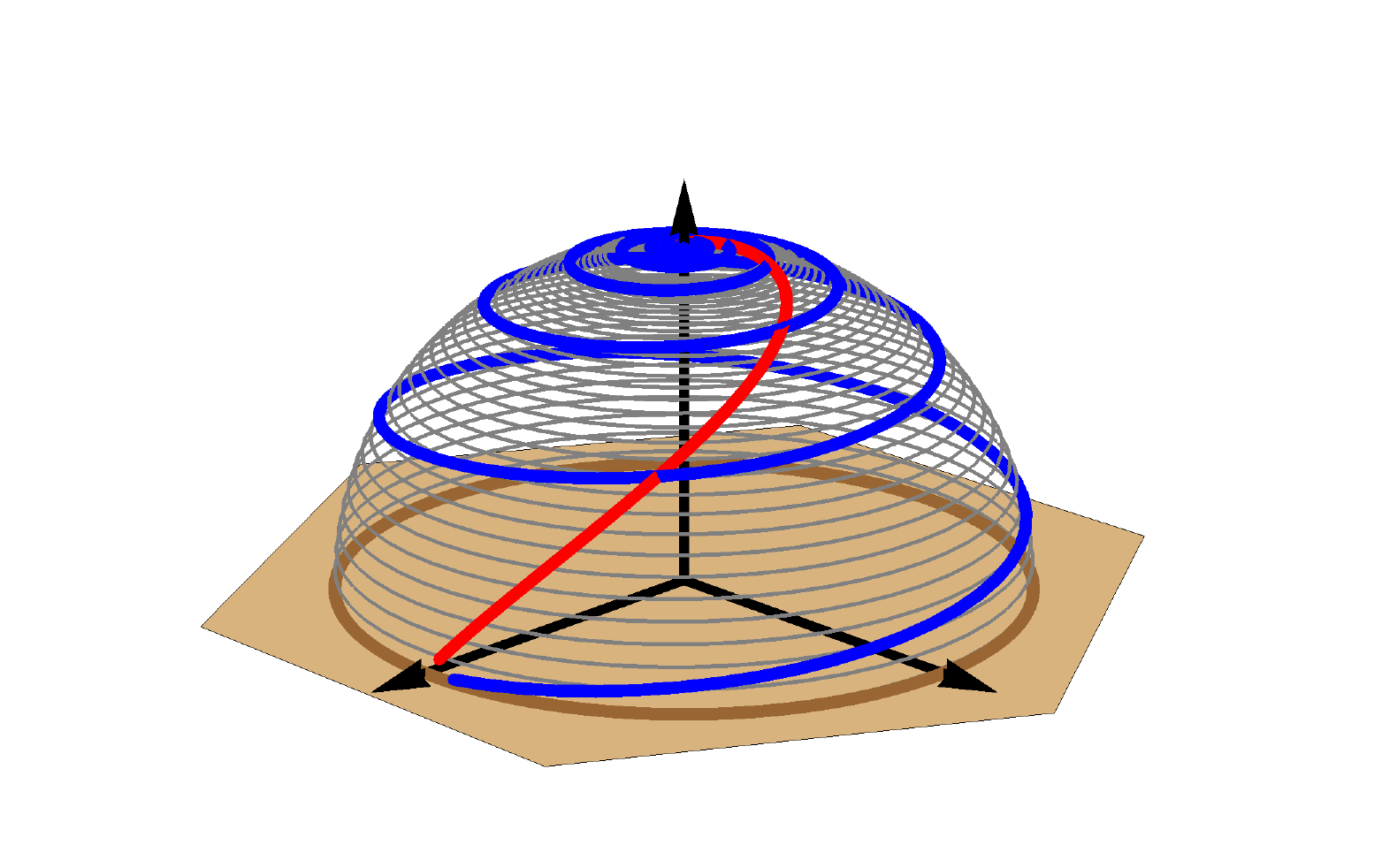}}
\put( 0.1,  2.7){\includegraphics[width=0.12\textwidth]{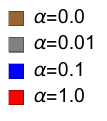}}
\put( 8.5,  0.0){\includegraphics[width=0.45\textwidth]{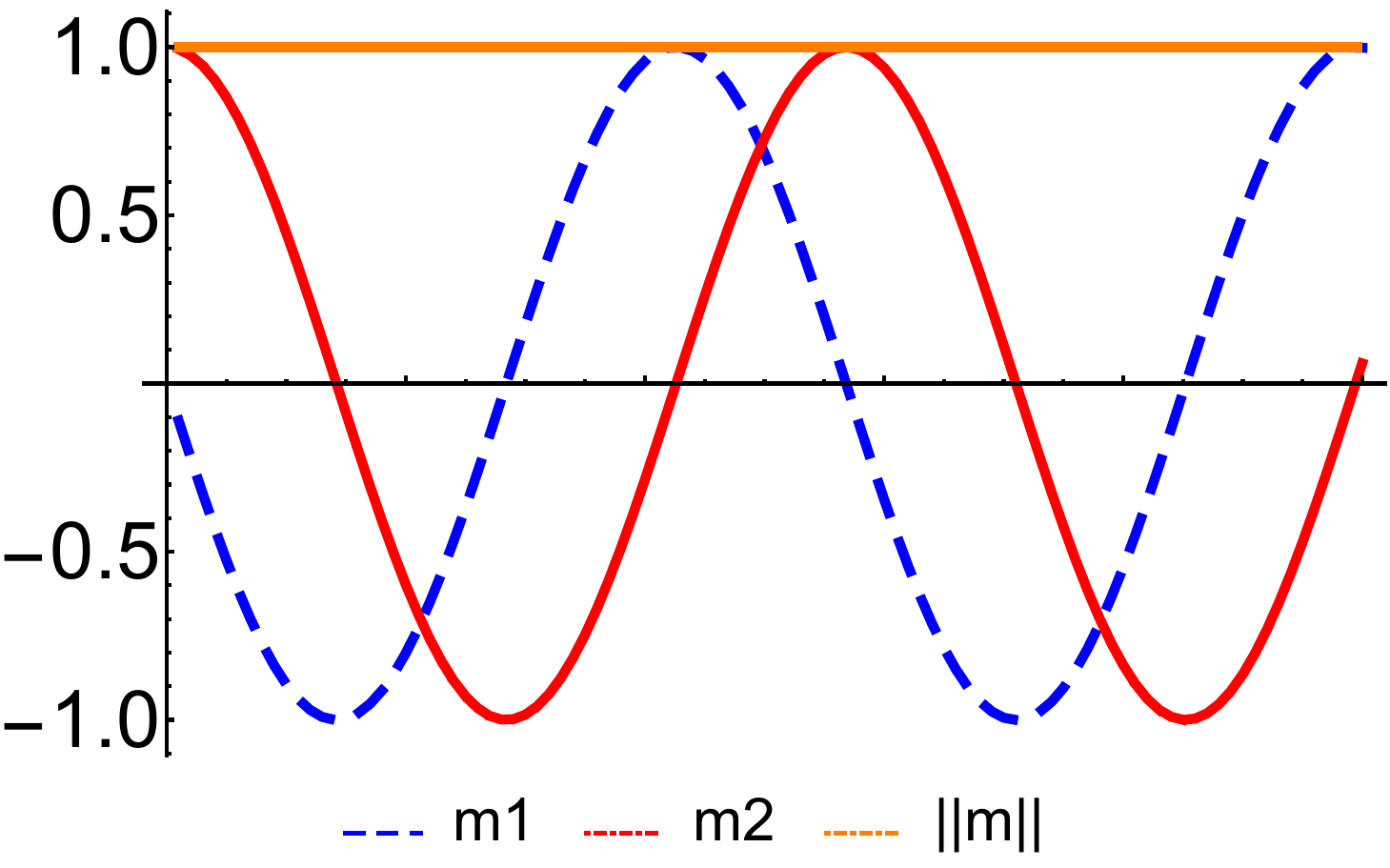}}
\put( 0.4, -0.5){a)}
\put( 8.2, -0.5){b)}
\put( 1.2,  0.3){$x_1$}
\put( 7.0,  0.5){$x_2$}
\put( 4.0,  4.8){$x_3$}
\put(15.5,  2.){t}
\end{picture}
\vspace*{5mm}
\caption{
The evolution path of the magnetization vectors is shown in a) for different damping coefficients, while b) shows the evolution of the magnetization components of the magnetization vector for a damping coefficient $\alpha=0$. The length of the vector remains constant $||\bbm||=1$ during the simulation.
} 
\label{fig:damped_switch}
\end{Figure}

\subsection{Stepsize stability} 
Stable algorithms should always guarantee the correct solution and thus physical correctness even with increasing time steps.
To check the properties of the presented algorithm we set the damping of the magnetization vectors to $\alpha=0$ in this example.
Hence, the externally applied magnetic field $\bH^{\rm ext}$ can  only cause, analogously to Fig.~\ref{fig:damped_switch}, a dissipation-free precession of the magnetization vectors. 
If a small time increment $\Delta t\approx0$ is used to simulate the precession and the position of the vector tip is stored after each converged time step, an almost perfect circle with diameter d=1 can be drawn. Thus $\norm{\bbm}=1$ applies correctly.
For increasingly larger time steps, the connecting line of the collected positions can deviate more and more from the geometry of the circle, since the magnetization also performs correspondingly larger rotations.
Regardless of this, however, the length of the magnetization vector must always satisfy $\norm{\bbm}=1$. This means, the tip of the vector must not leave the unit circle.  
This is analyzed in the following for the backward Euler method, the renormalized backward Euler method as well as for the presented exponential update algorithm. 
Under the externally applied field $\bH^{\rm ext}=10^{6}$A/m, the magnetization vector requires about t$_{\rm max}$~=~2.84$\times10^{-10}$s for a complete (360$^\circ$) rotation around the x$_3$ axis, starting from an initial configuration $\bbm=[1,0,0]^{\rm T}$. 
This rotation is performed in $k$ steps with $k=\{100,20,10,6,4,3,2\}$ in the study considered here.
This results in corresponding time increments of $\Delta t=t_{\rm max}/k$.
The starting point (initial configuration of magnetization) is highlighted with a red circle.
\begin{Figure}[h]
\begin{picture}(0,7.7)
\unitlength1cm
\put( 0.0,  7.5){\textbf{Backward Euler:}}
\put( 0.0,  6.8){2D projection}
\put( 7.0,  6.8){Norm $\bbm$}
\put(-0.5, -0.5){\rotatebox{0}{ \includegraphics[width=0.45\textwidth]{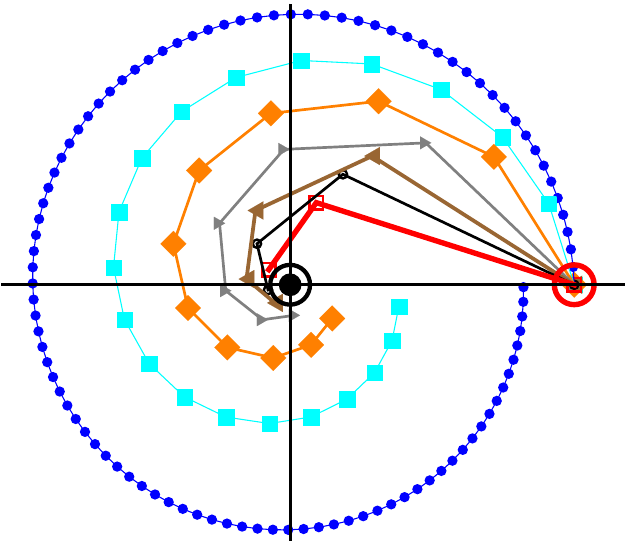}}}
\put( 7.0,  2.0){\includegraphics[width=0.55\textwidth]{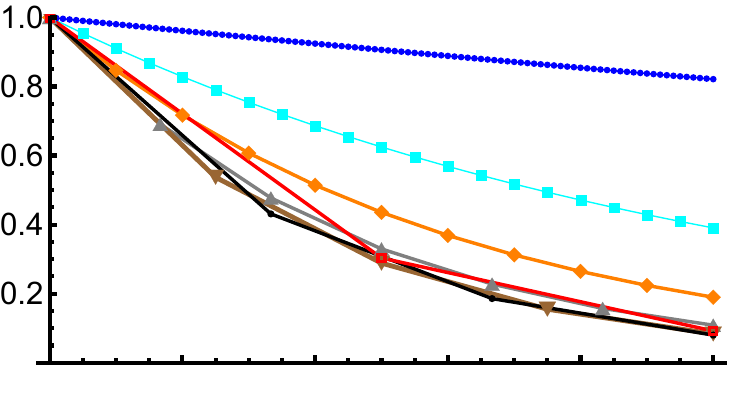}}
\put( 6.5, -1.5){\includegraphics[width=0.65\textwidth]{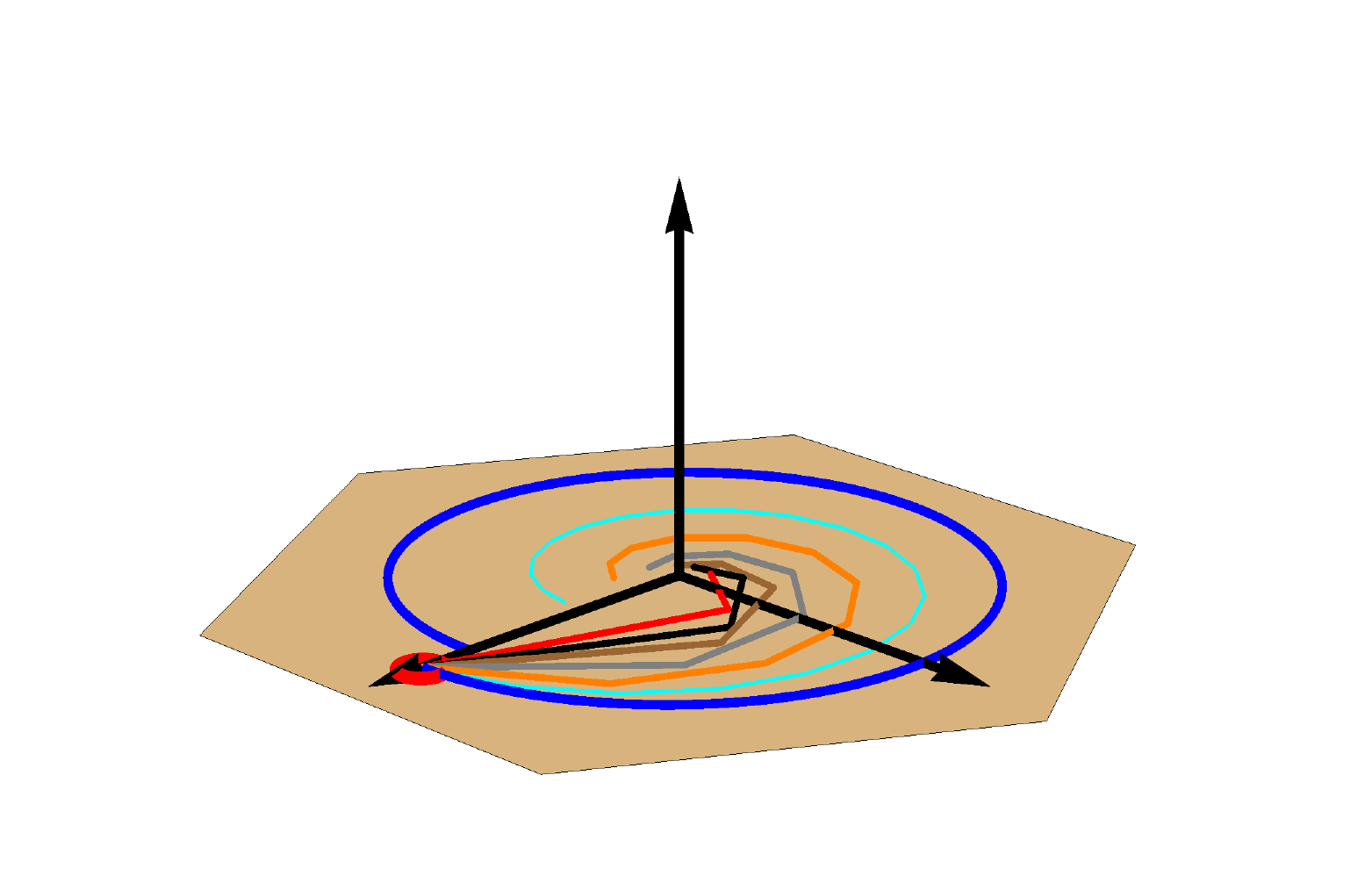}}
\put( 2.0, -2.2){\rotatebox{0}{ \includegraphics[width=0.7\textwidth]{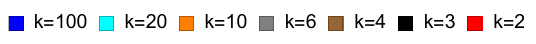}}}
\put( 0.4, -1.0){a)}
\put( 7.5,  1.5){b)}
\put( 8.2, -0.5){c)}
\put( 6.7,  2.5){$x_1$}
\put( 3.0,  6.5){$x_2$}
\put( 2.8,  3.3){$x_3$}
\put(15.5,  2.){t}
\put( 9.0, -0.5){$x_1$}
\put(14.0, -0.2){$x_2$}
\put(11.8,  3.0){$x_3$}
\end{picture}
\vspace*{20mm}
\caption{
a) and c) depict the precession of a magnetization vector within the x$_1$-x$_2$-plane for different time increments. The corresponding length of the magnetization vector during the rotation is presented in b).
} 
\label{fig:undamped_switch_euler}
\end{Figure}

\begin{Figure}[h]
\begin{picture}(0,8.0)
\unitlength1cm
\put( 0.0,  7.5){\textbf{Renormalized Backward Euler:}}
\put( 0.0,  6.8){2D projection}
\put( 7.0,  6.8){Norm $\bbm$}
\put(-0.5, -0.5){\rotatebox{0}{ \includegraphics[width=0.45\textwidth]{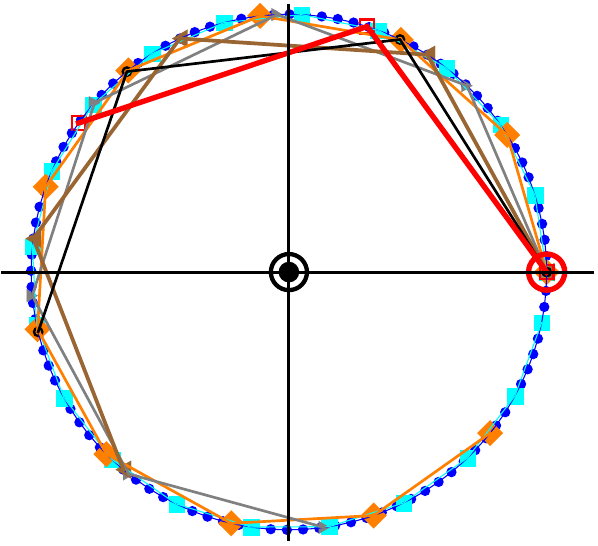}}}
\put( 7.0,  2.0){\includegraphics[width=0.55\textwidth]{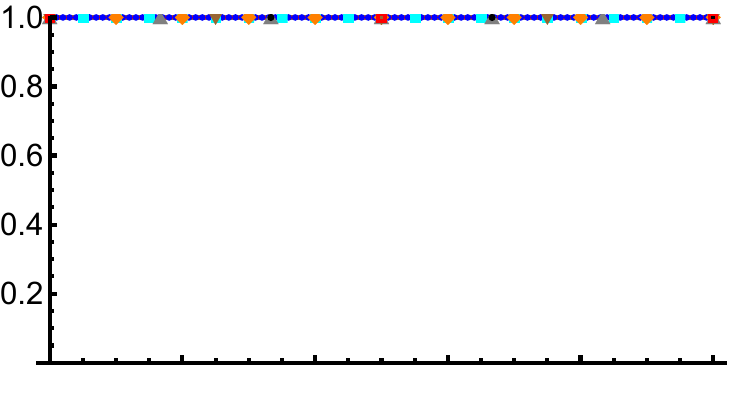}}
\put( 6.5, -1.5){\includegraphics[width=0.65\textwidth]{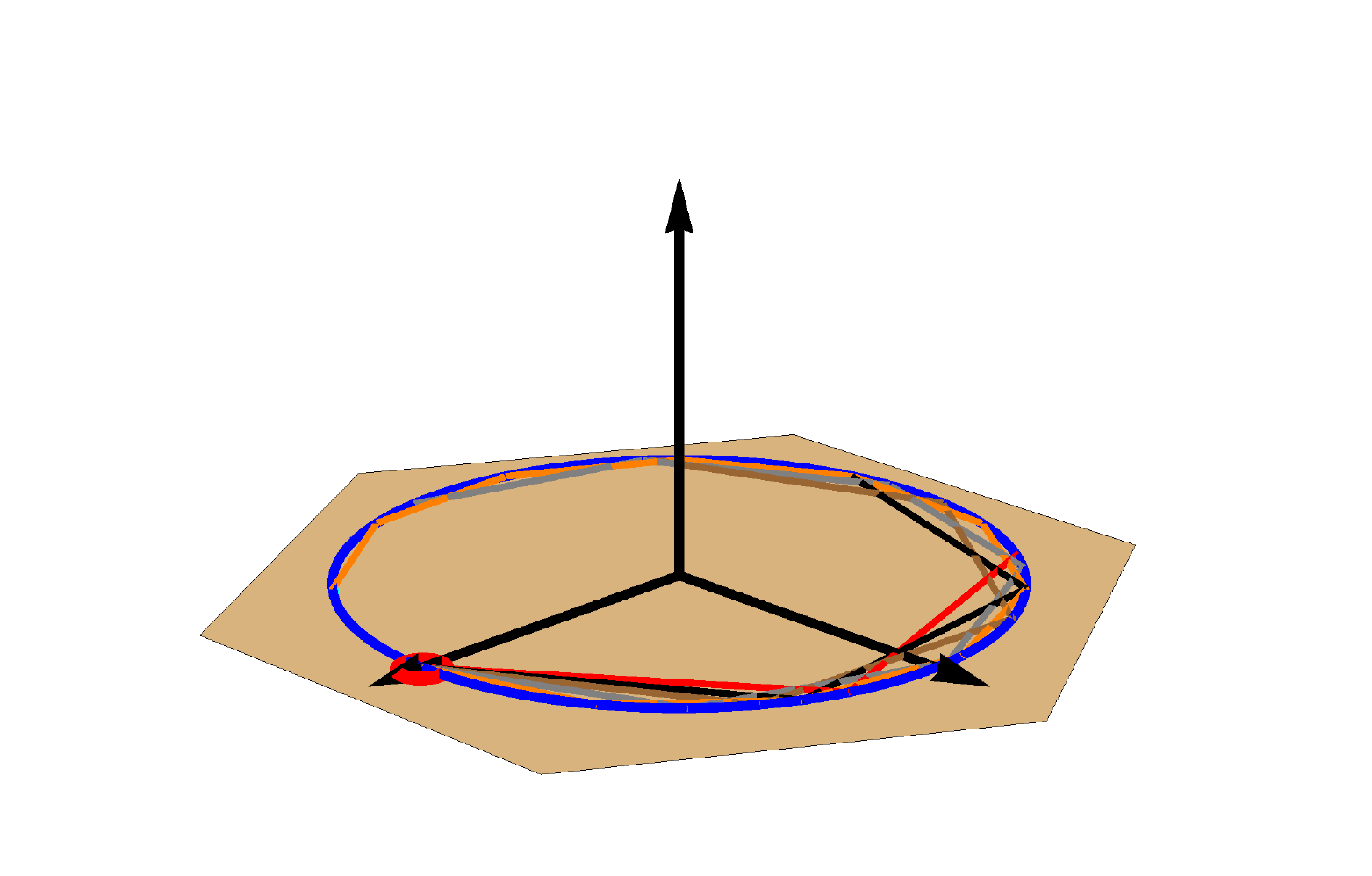}}
\put( 2.0, -2.2){\rotatebox{0}{ \includegraphics[width=0.7\textwidth]{figures/time_legend.pdf}}}
\put( 0.4, -1.0){a)}
\put( 7.5,  1.5){b)}
\put( 8.2, -0.5){c)}
\put( 6.7,  2.5){$x_1$}
\put( 3.0,  6.5){$x_2$}
\put( 2.8,  3.3){$x_3$}
\put(15.5,  2.){t}
\put( 9.0, -0.5){$x_1$}
\put(14.0, -0.2){$x_2$}
\put(11.8,  3.0){$x_3$}
\end{picture}
\vspace*{20mm}
\caption{
a) and c) depict the precession of a magnetization vector within the x$_1$-x$_2$-plane for different time increments.  The miss-match in the start and finish locations of the magnetization vectors can be clearly seen. The preservation of the norm during the rotation of the vector is presented in b).
} 
\label{fig:undamped_switch_euler_renorm}
\end{Figure}

The behavior of the magnetization vectors of the three methods is shown separately in Fig.~\ref{fig:undamped_switch_euler}, Fig.~\ref{fig:undamped_switch_euler_renorm} and Fig.~\ref{fig:undamped_switch_exp} as rotation around the x$_3$ axis.
Additionally, the norm $\norm{\bbm}$ is shown as a function of the simulated time for each applied number of time increments, to detect deviations form the unit length of the magnetization vectors. 
%


\begin{Figure}[h]
\begin{picture}(0,8.3)
\unitlength1cm
\put( 0.0,  7.5){\textbf{Implicite Midpoint Rule:}}
\put( 0.0,  6.8){2D projection}
\put( 7.0,  6.8){Norm $\bbm$}
\put(-0.5, -0.5){\rotatebox{0}{ \includegraphics[width=0.45\textwidth]{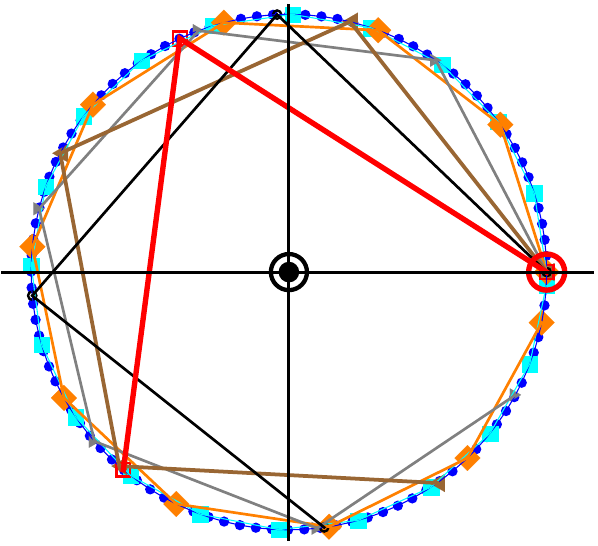}}}
\put( 7.0,  2.0){\includegraphics[width=0.55\textwidth]{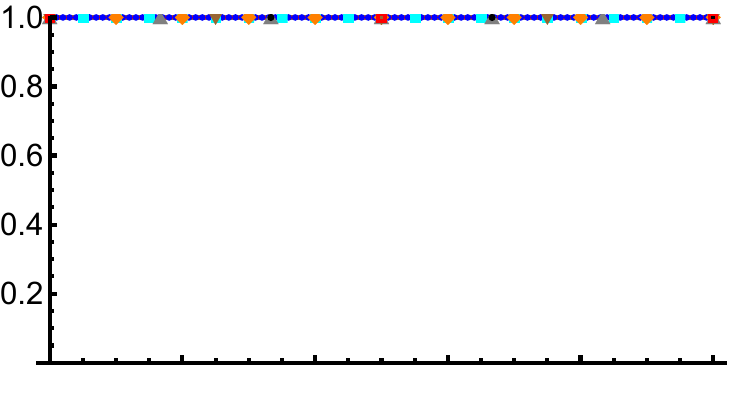}}
\put( 6.5, -1.5){\includegraphics[width=0.65\textwidth]{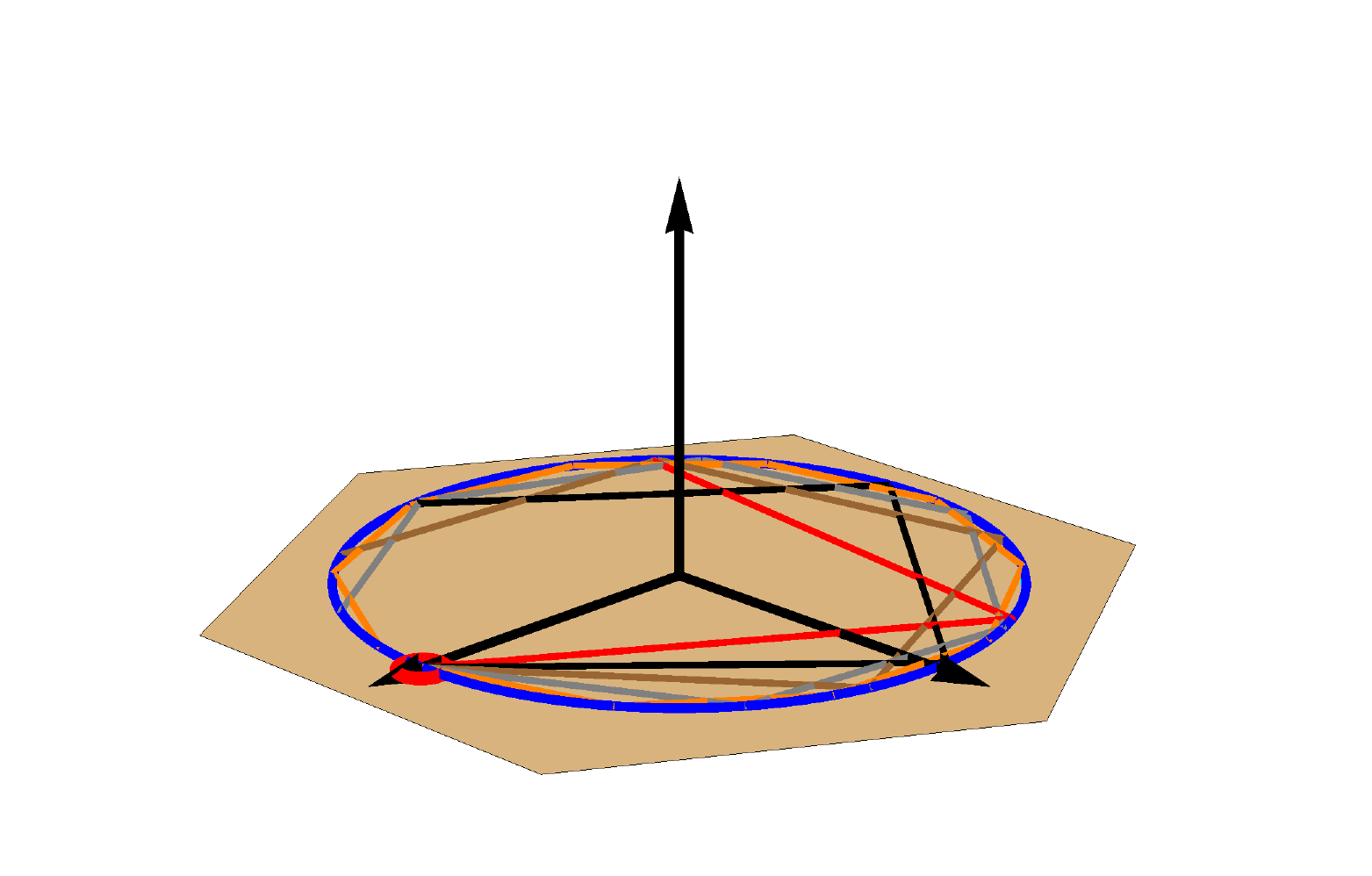}}
\put( 2.0, -1.9){\rotatebox{0}{ \includegraphics[width=0.6\textwidth]{figures/time_legend.pdf}}}
\put( 0.4, -1.0){a)}
\put( 7.5,  1.5){b)}
\put( 8.2, -0.5){c)}
\put( 6.7,  2.5){$x_1$}
\put( 3.0,  6.2){$x_2$}
\put( 2.6,  3.2){$x_3$}
\put(15.5,  2.){t}
\put( 9.0, -0.5){$x_1$}
\put(14.0, -0.2){$x_2$}
\put(11.8,  3.0){$x_3$}
\end{picture}
\vspace*{15mm}
\caption{
a) depicts the precession of a magnetization vector within the x$_1$-x$_2$-plane for different time increments. The corresponding norm preservation during the rotation of the vector is presented in b).
} 
\label{fig:undamped_switch_MPR}
\end{Figure}

As can be clearly seen in Fig.~\ref{fig:undamped_switch_euler}, the implicit Euler method is not length preserving. This means that the length of the magnetization vector, as rotation and time proceeds, deviates from the unit length and becomes smaller. 
The difference between the initial length and that of the last point in time is also strongly dependent on the size of the time increments. 
Thus, the smallest deviation of the length from one is observed at $k=100$ with $\norm{\bbm}=0.82$ and the strongest deviation at $k=3$ with $\norm{\bbm}=0.08$.  
This problem can be solved for the considered example by renormalizing the vectors after each completed time step. Corresponding results presented in Fig.~\ref{fig:undamped_switch_euler_renorm} indicate the preservation of the unit length constraint of the magnetization vector. 
However, in \cite{ReiXuSch:2022:acs} the authors show that the renormalized residuals do not match the converged residuals leading to the quintessence, that the renormalized solution is far from the minimized/calculated solution.
A much better choice of parameterization at this point is the exponential update algorithm. 
As can be seen in Fig.~\ref{fig:undamped_switch_exp}b, this algorithm represents a length-preserving and thus rotation-exact method. 
%
Besides the missing ability to keep the vector length constant the backward Euler methods also suffers from deviating positions of the vectors compared to reference solutions. 
While the starting position of the vectors should also correspond to the final position both, backwards Euler methods fail to do so. 
This deviation from the correct position can already be seen in Fig.~\ref{fig:undamped_switch_euler} and Fig.~\ref{fig:undamped_switch_euler_renorm}, but Fig.~\ref{fig:StepSizePositionError} indicates a direct correlation of the position of the vector to the applied step size during the simulation. 
\begin{Figure}[h]
\begin{picture}(0,7.5)
\unitlength1cm
\put( 0.0,  7.5){\textbf{Exponential update:}}
\put( 0.0,  6.8){2D projection}
\put( 7.0,  6.8){Norm $\bbm$}
\put(-0.5, -0.5){\rotatebox{0}{ \includegraphics[width=0.45\textwidth]{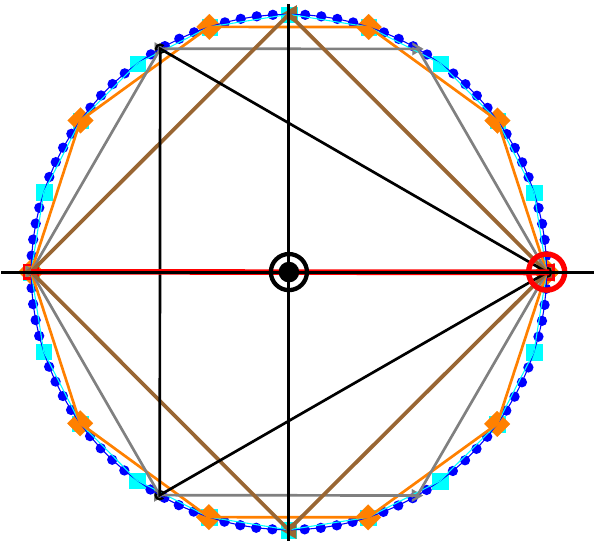}}}
\put( 7.0,  2.0){\includegraphics[width=0.55\textwidth]{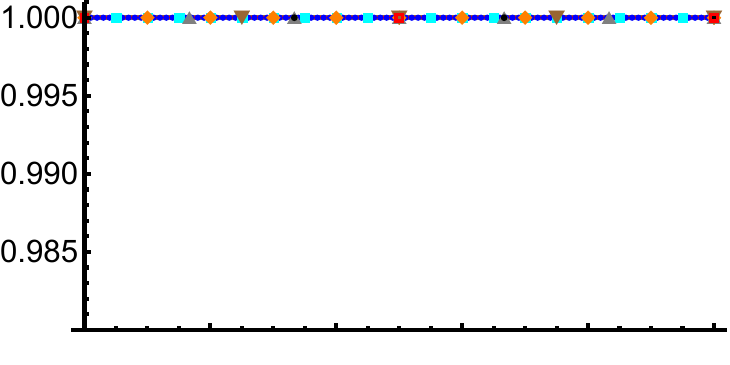}}
\put( 6.5, -1.5){\includegraphics[width=0.65\textwidth]{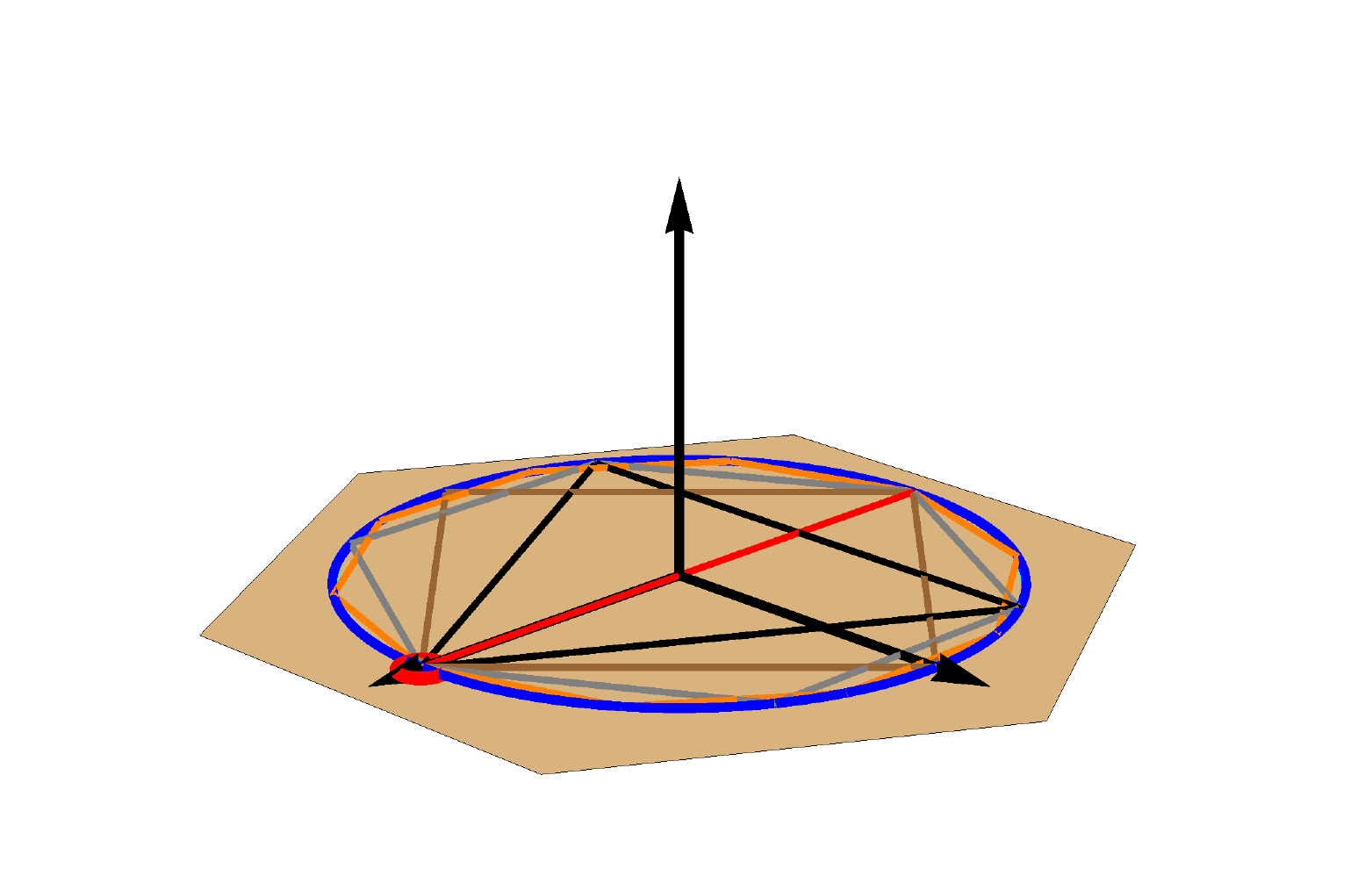}}
\put( 2.0, -1.9){\rotatebox{0}{ \includegraphics[width=0.6\textwidth]{figures/time_legend.pdf}}}
\put( 0.4, -1.0){a)}
\put( 7.5,  1.5){b)}
\put( 8.2, -0.5){c)}
\put( 6.7,  2.5){$x_1$}
\put( 3.0,  6.2){$x_2$}
\put( 2.6,  3.2){$x_3$}
\put(15.5,  2.){t}
\put( 9.0, -0.5){$x_1$}
\put(14.0, -0.2){$x_2$}
\put(11.8,  3.0){$x_3$}
\end{picture}
\vspace*{15mm}
\caption{
a) depicts the precession of a magnetization vector within the x$_1$-x$_2$-plane for different time increments. The corresponding norm preservation during the rotation of the vector is presented in b).
} 
\label{fig:undamped_switch_exp}
\end{Figure}

\vspace{-1mm}
Larger time increments fail to reach the initial position within the considered time t$_{\rm max}$, while smaller time steps $\Delta t \approx 0$ lead to better results but also increasing computational times.
This behavior contrasts the exponential update algorithm.
As depicted in Fig.~\ref{fig:undamped_switch_exp} and Fig.~\ref{fig:StepSizePositionError} the tip of the magnetization vector never leaves the unit circle, not even for the coarsest time discretization, and always ends at the intended point.
As a comparison of Fig.~\ref{fig:undamped_switch_euler}, Fig.~\ref{fig:undamped_switch_euler_renorm}, and Fig.~\ref{fig:StepSizePositionError} shows, the exponential delivers the most accurate solution. 
Hence $\bbm^*$ indicates the last position of the magnetization vector obtained from a simulation using the exponential update algorithm with an increment step size of $\Delta=t_{\rm max}/1000$.

%
\vspace{50mm}
\begin{minipage}{0.49\textwidth}
   \begin{picture}(6.0,0.0)
\unitlength1cm
\put( 0.0, 0.5){\includegraphics[width=7.5cm]{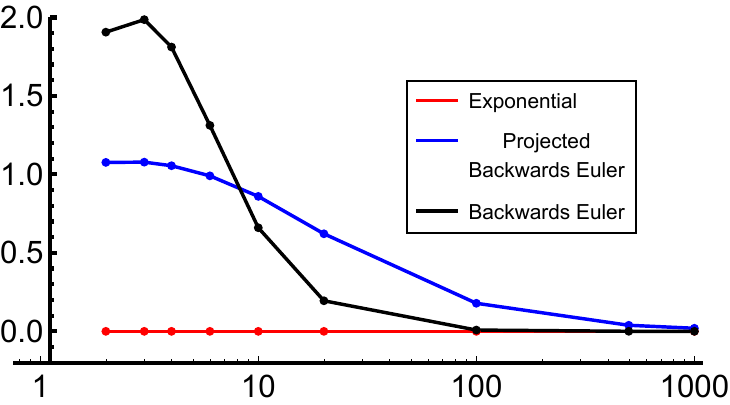}}
\put(-0.5, 3.0){$\varphi$}
\put( 1.5, 0.0){current number of time step}
\end{picture}
    \captionof{figure}{The deviation of the position of the magnetization vector $\bbm$ is shown with respect to a reference solution $\bbm^*$ at time t$_{\rm max}$~=~2.84$\times10^{-10}$s. The reference solution was generated with the exponential update and a number of 1000 time increments.}
    \label{fig:StepSizePositionError}
\end{minipage}
\hspace{2mm}
\begin{minipage}{0.49\textwidth}
    \centering
   \begin{picture}(7.5,0.0)
\unitlength1cm
\put( 0.0, 1.5){\includegraphics[width=7.5cm]{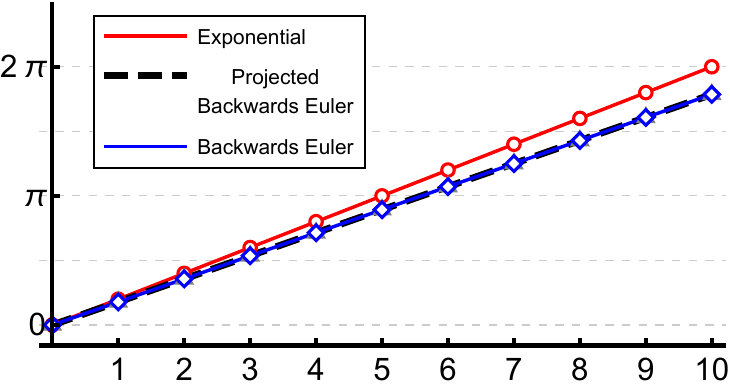}}
\put(-0.2, 3.5){$\varphi$}
\put( 1.5, 0.95){current number of time step}
\end{picture}
\vspace{-9mm}
    \captionof{figure}{Evolution of the polar angle $\varphi=\tan^{-1}(m_2/m_1)$ of the backwards Euler, the renormalized backwards Euler, and the exponential update scheme.}
    \label{fig:StepSizeAngleError}
\end{minipage}
%
%
The deviation of the magnetization vector from its correct position can be described by the angle $\varphi$ between the x$_1$-axis and the vector itself. 
This deviation is shown for the case of $\Delta=t_{\rm max}/10$ in Fig.~\ref{fig:StepSizeAngleError}. 
Here it becomes quite evident that the exponential update algorithm finally returns to its initial position, while the backward Euler algorithms, i.e., the normal as well as the renormalized one, does not reach this position. 
Likewise it is to be recognized in Fig.~\ref{fig:StepSizeAngleError} that the length of the vector has no influence on the rotation angle. 
%
%

\begin{Figure}[h]
\begin{picture}(0, 12.)
\unitlength1cm
\put( 0.2,  9.5){\includegraphics[width=0.45\textwidth]{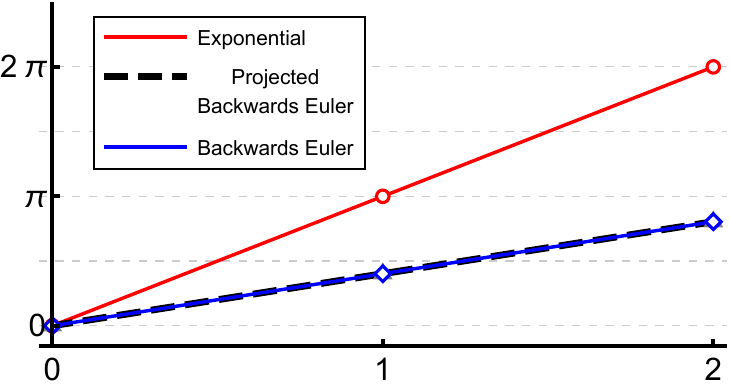}}
\put( 8.5,  9.5){\includegraphics[width=0.45\textwidth]{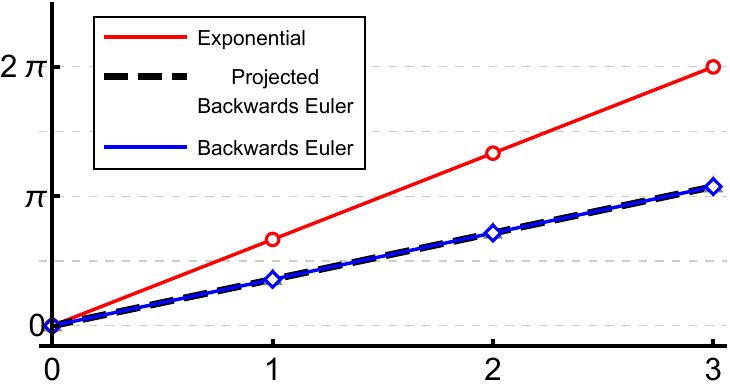}}
\put( 0.2,  4.3){\includegraphics[width=0.45\textwidth]{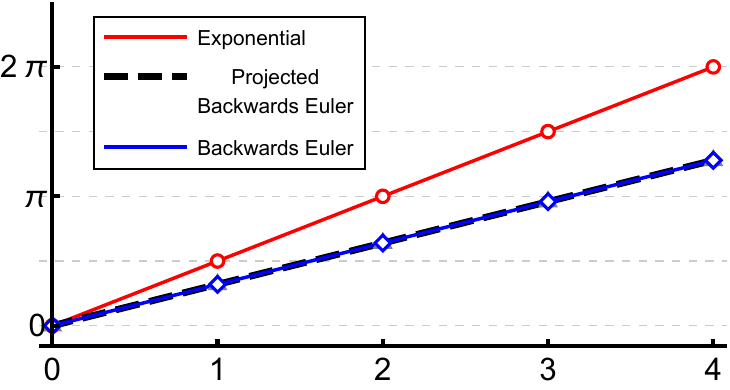}}
\put( 8.5,  4.3){\includegraphics[width=0.45\textwidth]{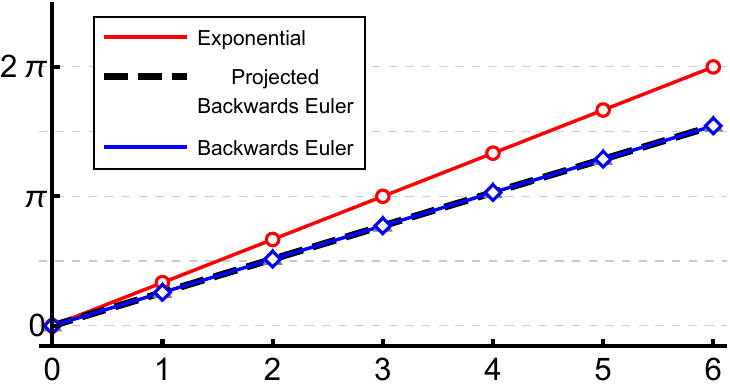}}
\put( 0.2, -0.9){\includegraphics[width=0.45\textwidth]{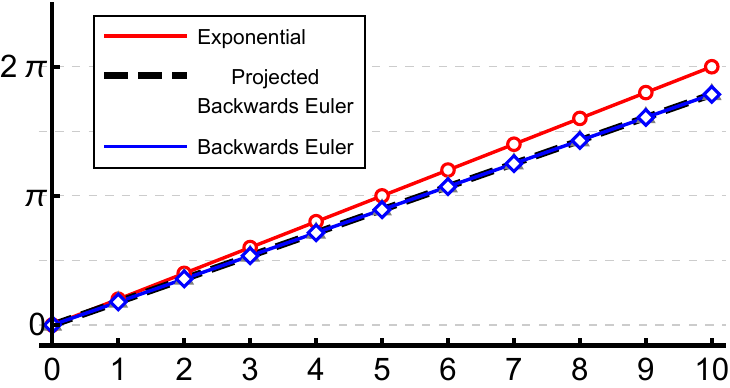}}
\put( 8.5, -0.9){\includegraphics[width=0.45\textwidth]{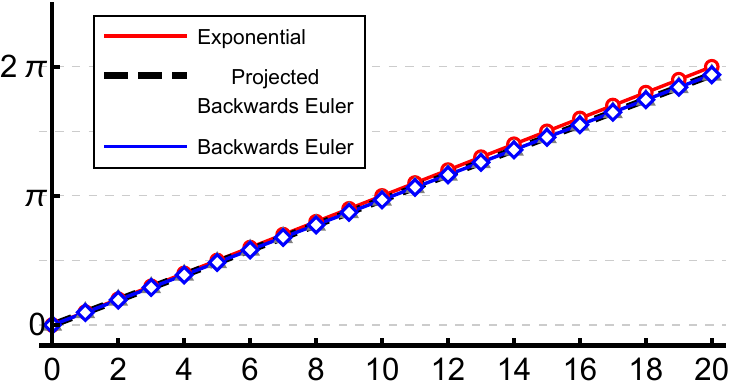}}
\put( 0.0, 11.4){$\varphi$}
\put( 4.8, 10.3){k=2}
\put( 1.8,  9.0){current number of time step k}
\put( 8.2, 11.4){$\varphi$}
\put(13.0, 10.3){k=3}
\put(10.0,  9.0){current number of time step k}
\put( 0.0,  6.2){$\varphi$}
\put( 4.8,  5.3){k=4}
\put( 1.8,  3.8){current number of time step k}
\put( 8.2,  6.3){$\varphi$}
\put(13.0,  5.2){k=6}
\put(10.0,  3.8){current number of time step k}
\put( 0.0,  1.0){$\varphi$}
\put( 4.8,  0.1){k=10}
\put( 1.8, -1.4){current number of time step k}
\put( 8.2,  1.0){$\varphi$}
\put(13.0,  0.1){k=20}
\put(10.0, -1.4){current number of time step k}
\end{picture}
\vspace*{20mm}
\caption{
Evolution of the polar angle $\varphi=\tan^{-1}(m_2/m_1)$ of the backwards Euler, the renormalized backwards Euler, and the exponential update scheme. 
} 
\label{fig:StepSizeAngleErrorManyDt}
\end{Figure}

Fig.~\ref{fig:StepSizeFinalPositionError} shows the norm of the last calculated magnetization vector for all three considered algorithms, plotted over the number of iterations k used. 
The last position should be identical to the initial position. 

\begin{Figure}[h]
\begin{picture}(0,3.0)
\unitlength1cm
\put( 4.0, -1.5){\includegraphics[width=0.5\textwidth]{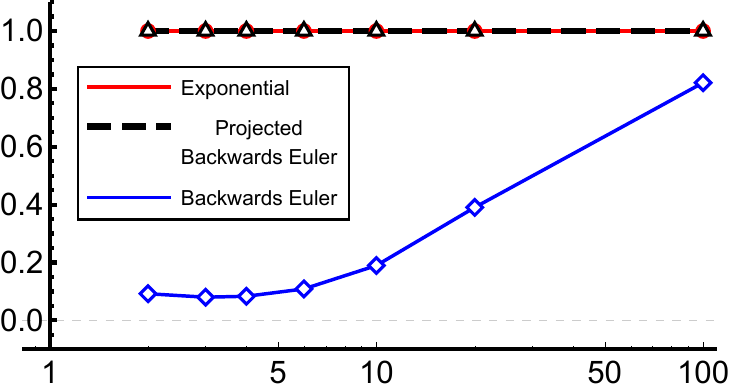}}
\put( 2.5,  0.8){$\norm{\bbm}$}
\put( 5.5, -2.0){number of total time increments k}
\end{picture}
\vspace*{20mm}
\caption{
Deviation of the unit length of the last calculated magnetization vector $\bbm$, at an angle of $2\pi$, for the exponential update, the renormalized backward Euler, and the classical backward Euler. 
} 
\label{fig:StepSizeFinalPositionError}
\end{Figure}

\pagebreak
\section{\hspace{-5mm}. Conclusion}
\vspace{-4mm}
In this work, the exponential update algorithm for the description of exact rotations is discussed in comparison to the backward Euler and the renormalized backward Euler method.
Micromagnetism with the Landau-Lifshitz-Gilbert equation represents the necessary physical frame work, while
the characteristic feature of this theory is the magnetization vector that must always satisfy the non-convex condition $\norm{\bbm}=1$.
While both the classical and the renormalized backward Euler give bad results for increased time increments, regarding the conservation of length as well as for the position of the vector on the unit sphere, the exponential update algorithm proves to be consistently stable.
The advantages of the exponential update algorithm are therefore obvious. 
Larger time steps with the same level of quality lead to significantly shorter computing times.

{\bf Acknowledgements} We gratefully acknowledge the financial support of the German Research Foundation
(DFG) in the framework of the CRC/TRR 270, project A07 "Scale-bridging of magneto-mechanical mesostructures of
additive manufactured and severe plastically deformed materials", project number 405553726. 

\bibliographystyle{plainnat}
\bibliography{magneto_17}

@Article{BarPro:2007:cpi,
  author  = {S. Bartels and A. Prohl},
  title   = {Constraint preserving implicit finite element discretization of harmonic map flow into spheres},
  journal = {Mathematics of Computation},
  year    = {2007},
  volume  = {76},
  number  = {260},
  pages   = {1847-1859},
}

@PhdThesis{Eid:2004:aim,
  author = {B. Eidel},
  title  = {Anisotropic Inelasticity - Modelling, Simulation, Validation},
  school = {Technical Universit\"at Darmstadt},
  year   = {2004},
}

@Book{HaiLabWan:2006:gni,
  title     = {IGeometric Numerical Integration, Structure-Preserving Algorithms for Ordinary Differential Equations},
  publisher = {Springer},
  year      = {2006},
  author    = {E.~Hairer and  C.~Lubich and G.~Wanner},
  edition   = {2},
}

@Book{CulGra:2009:itm,
  title     = {Introduction to magnetic materials},
  publisher = {John Wiley and Sons},
  year      = {2009},
  author    = {B.D.~Cullity and C.D.~Graham},
  edition   = {2},
}

@Article{YiXu:2014:acf,
  author  = {M. Yi and B.-X. Xu},
  title   = {A constraint-free phase field model for ferromagnetic domain evolution},
  journal = {Proceedings of the Royal Society A},
  year    = {2014},
  volume  = {470},
  number  = {2171},
  pages   = {20140517},
  doi     = {https://doi.org/10.1098/rspa.2014.0517},
}

@Article{FenVis:2002:saa,
  author  = {X.~Feng and P.B.~Visscher},
  journal = {Journal of Applied Physics},
  title   = {Stability and accuracy of {E}uler and quaternion micromagnetic algorithms},
  year    = {2002},
  number  = {10},
  pages   = {8712-8714},
  volume  = {91},
  doi     = {10.1063/1.1452284},
}

@Article{LewNig:2003:gio,
  author  = {D.~Lewis and N.~Nigam},
  title   = {Geometric integration on spheres and some interesting applications},
  journal = {Journal of Computational and Applied Mathematics},
  year    = {2003},
  volume  = {151},
  number  = {1},
  pages   = {141-170},
  doi     = {https://doi.org/10.1016/S0377-0427(02)00743-4},
}

@InCollection{LanLif:1935:ott,
  author    = {L. Landau and E. Lifshitz},
  title     = {On the theory of the dispersion of magnetic permeability in ferromagnetic bodies},
  booktitle = {Perspectives in Theoretical Physics},
  publisher = {Elsevier},
  year      = {1935},
  pages     = {51--65},
  doi       = {10.1016/B978-0-08-036364-6.50008-9},
}

@Book{Bro:1963:m,
  title     = {Micromagnetics},
  publisher = {Wiley},
  year      = {1963},
  author    = {W.F. Brown},
}

@Article{Gil:2004:apt,
  author    = {T.L. Gilbert},
  title     = {A phenomenological theory of damping in ferromagnetic materials},
  journal   = {IEEE transactions on magnetics},
  year      = {2004},
  volume    = {40},
  number    = {6},
  pages     = {3443--3449},
  publisher = {IEEE},
}

@Article{ReiXuSch:2022:acs,
  author  = {M.~Reichel and B.-X.~Xu and J.~Schr\"oder},
  title   = {A comparative study of finite element schemes for micromagnetic mechanically coupled simulations},
  journal = {Journal of Applied Physics},
  year    = {2022},
  volume  = {132},
  number  = {18},
  pages   = {183903},
  doi     = {10.1063/5.0105613},
}

@Book{Pro:2001:cm,
  title     = {Computational micromagnetism},
  publisher = {Springer},
  year      = {2001},
  author    = {A.~Prohl},
  edition   = {1},
  doi       = {10.1007/978-3-663-09498-2},
}

@Article{VisFen:2002:qba,
  author  = {P.B.~Visscher and X.~Feng},
  title   = {Quaternion-based algorithm for micromagnetics},
  journal = {Physical Review B},
  year    = {2002},
  volume  = {65},
  number  = {10},
  pages   = {104412},
  doi     = {10.1103/PhysRevB.65.104412},
}

@Article{KruPro:2006:rdi,
  author  = {M.~Kru{\v{z}}{\'{i}}k and A.~Prohl},
  title   = {Recent Developments in the Modeling, Analysis, and Numerics of Ferromagnetism},
  journal = {SIAM Review},
  year    = {2006},
  volume  = {48},
  number  = {3},
  pages   = {439--483},
  doi     = {10.1137/S0036144504446187},
}

@Article{SriKeiMie:2015:chi,
  author    = {A. Sridhar and M.A. Keip and C. Miehe},
  title     = {Computational Homogenization in Micro-Magneto-Elasticity},
  journal   = {Proceedings in Applied Mathematics and Mechanics},
  year      = {2015},
  volume    = {15},
  number    = {1},
  pages     = {363--364},
  doi       = {10.1002/pamm.201510172},
  publisher = {Wiley Online Library},
}

@Article{Lan:2008:act,
  author    = {C.M. Landis},
  title     = {A continuum thermodynamics formulation for micro-magneto-mechanics with applications to ferromagnetic shape memory alloys},
  journal   = {Journal of the Mechanics and Physics of Solids},
  year      = {2008},
  volume    = {56},
  number    = {10},
  pages     = {3059--3076},
  doi       = {10.1016/j.jmps.2008.05.004},
  publisher = {Elsevier},
}

@Article{WanZah:2013:ars,
  author  = {J. Wang and J. Zhang},
  title   = {A real-space phase field model for the domain evolution of ferromagnetic materials},
  journal = {International Journal of Solids and Structures},
  year    = {2013},
  volume  = {50},
  number  = {22-23},
  pages   = {3597-3609},
  doi     = {doi.org/10.1016/j.ijsolstr.2013.07.001},
}

@Article{ZhaZhaPei:2016:afe,
  author    = {H.~Zhang and X.~Zhang and Y.~Pei},
  title     = {A finite element based real-space phase field model for domain evolution of ferromagnetic materials},
  journal   = {Computational Materials Science},
  year      = {2016},
  volume    = {118},
  pages     = {214--223},
  doi       = {10.1016/j.commatsci.2016.03.020},
  publisher = {Elsevier},
}

@Article{OhmYiGutXu:2022:pfm,
  author    = {D.~Ohmer and M.~Yi and O.~Gutfleisch and B.~X.~Xu},
  title     = {Phase-field modelling of paramagnetic austenite--ferromagnetic martensite transformation coupled with mechanics and micromagnetics},
  journal   = {International Journal of Solids and Structures},
  year      = {2022},
  volume    = {238},
  pages     = {111365},
  doi       = {10.1016/j.ijsolstr.2021.111365},
  publisher = {Elsevier},
}

@Article{SzaBudTouFru:acf:2008:acf,
  author  = {H. Szambolics and L.D. Buda-Prejbeanu and J.C. Toussaint and O. Fruchart},
  journal = {Computational material science},
  title   = {A constrained finite element formulation for the {L}andau-{L}ifshitz-{G}ilbert equations},
  year    = {2008},
  number  = {2},
  pages   = {253--258},
  volume  = {44},
  doi     = {10.1016/j.commatsci.2008.03.01},
}

@Article{FidSch:2000:mmt,
  author  = {J. Fidler and T. Schrefl},
  title   = {Micromagnetic modelling-the current state of the art},
  journal = {Journal of Physics D: Applied Physics},
  year    = {2000},
  volume  = {33},
  pages   = {R135-R156},
}

@Article{SueSchFid:2000:mso,
  author    = {D.~S\"uss and T.~Schrefl and J.~Fidler},
  title     = {Micromagnetics simulation of high energy density permanent magnets},
  journal   = {IEEE transactions on magnetics},
  year      = {2000},
  volume    = {36},
  number    = {5},
  pages     = {3282--3284},
  doi       = {10.1109/20.908770},
  publisher = {IEEE},
}

@Article{SchFidSchSueForTsi:2003:spm,
  author    = {W. Scholz and J. Fidler and T. Schrefl and D. S\"uss and H. Forster and V. Tsiantos},
  title     = {Scalable parallel micromagnetic solvers for magnetic nanostructures},
  journal   = {Computational Materials Science},
  year      = {2003},
  volume    = {28},
  number    = {2},
  pages     = {366--383},
  doi       = {10.1016/S0927-0256(03)00119-8},
  publisher = {Elsevier},
}

@Article{DorSchXuKeiMue:2018:cpf,
  author  = {W. Dornisch and D. Schrade and B.-X. Xu and M.-A. Keip and R. M\"uller},
  title   = {Coupled phase field simulations of ferroelectric and ferromagnetic layers in multiferroic heterostructures},
  journal = {Archive of Applied Mechanics},
  year    = {2018},
  volume  = {89},
  number  = {6},
  pages   = {1031-1056},
  doi     = {https://doi.org/10.1007/s00419-018-1480-9},
}

@Article{MieEth:2014:agc,
  author  = {C. Miehe and G. Ethiraj},
  title   = {A geometrically consistent incremental variational formulation for phase field models in micromagnetics},
  journal = {Computer Methods in Applied Mechanics and Engineering},
  year    = {2012},
  volume  = {245-246},
  pages   = {331-347},
  doi     = {https://doi.org/10.1016/j.cma.2012.03.021},
}

@Article{Alo:2008:anf,
  author    = {F.~Alouges},
  title     = {A new finite element scheme for {L}andau-{L}ifchitz equations},
  journal   = {Discrete \& Continuous Dynamical Systems-S},
  year      = {2008},
  volume    = {1},
  number    = {2},
  pages     = {187},
  publisher = {American Institute of Mathematical Sciences},
}

@Article{dAqSerMia:2005:gio,
  author  = {M.~d'Aquino and C.~Serpico and G.~Miano},
  title   = {Geometrical integration of {L}andau-{L}ifshitz-{G}ilbert equation based on the mid-point rule},
  journal = {Journal of Computational Physics},
  year    = {2005},
  volume  = {209},
  number  = {2},
  pages   = {730-753},
  doi     = {10.1016/j.jcp.2005.04.001},
}

@Article{Cim:2008:tll,
  author  = {I. Cimr{\'a}k},
  journal = {Archives of Computational Methods in Engineering},
  title   = {A Survey on the numerics and computations for the {L}andau-{L}ifshitz equation of micromagnetism},
  year    = {2008},
  number  = {3},
  pages   = {1--37},
  volume  = {15},
  doi     = {10.1007/BF03024947},
}

@Article{Abe:2019:mas,
  author    = {C.~Abert},
  title     = {Micromagnetics and spintronics: models and numerical methods},
  journal   = {The European Physical Journal B},
  year      = {2019},
  volume    = {92},
  number    = {6},
  pages     = {1--45},
  doi       = {doi.org/10.1140/epjb/e2019-90599-6},
  publisher = {Springer},
}

@Book{Ber:1998:him,
  title     = {Hysteresis in Magnetism For Physicists, Materials Scientists, and Engineers},
  publisher = {Gulf Professional Publishing},
  year      = {1998},
  author    = {G. Bertotti},
  address   = {San Diego},
  doi       = {10.1016/B978-0-12-093270-2.X5048-X},
}

@Book{Coe:2010:mam,
  title     = {Magnetism and magnetic materials},
  publisher = {Cambridge university press},
  year      = {2010},
  author    = {J.M.D.~Coey},
  address   = {Cambridge},
  doi       = {10.1017/CBO9780511845000},
}

@Book{HubSch:2008:mdt,
  title     = {Magnetic domains: the analysis of magnetic microstructures},
  publisher = {Springer Science \& Business Media},
  year      = {2008},
  author    = {A. Hubert and R. Sch{\"a}fer},
  volume    = {3},
  address   = {Berlin},
}

@Book{KroPar:2007:hom,
  author    = {H. Kronm{\"u}ller and S.P. Parkin},
  publisher = {Wiley New York},
  title     = {Handbook of magnetism and advanced magnetic materials},
  year      = {2007},
  volume    = {2},
}

@article{DorWoeWul:2023:cbn,
    author = {C. Dorn and M. H\"orsting and S. Wulfinghoff},
    title = {Computing Barkhausen noise spectra for magnetostrictive thin film composites using efficient magnetization-magnitude preserving simulation techniques},
    journal = {Journal of Applied Physics},
    volume = {134},
    number = {13},
    pages = {133901},
    year = {2023},
    doi = {10.1063/5.0157906}
}

@Article{ReiNieSch:2023:ems,
  author  = {M. Reichel and R. Niekamp and J. Schr\"oder},
  journal = {Journal of Applied Physics},
  title   = {Efficient micromagnetic simulations based on a perturbed {L}agrangian function},
  year    = {2023},
  number  = {10},
  pages   = {103901},
  volume  = {134},
  doi     = {10.1063/5.0159273},
}


\end{document}